\documentclass{article}

\usepackage[english]{babel}
\usepackage{listings}
\makeatletter
\usepackage{jlcode}
\lstdefinestyle{code}{
  language=Julia, 
  showstringspaces=false,
  keywordstyle=\color{blue},
  commentstyle=\color{gray},
  identifierstyle=\color[RGB]{0,102,0},
  columns=fullflexible,
  keepspaces=true
}
\lstnewenvironment{code}{\lstset{style=code}}{}

\usepackage[letterpaper,top=2cm,bottom=2cm,left=3cm,right=3cm,marginparwidth=1.75cm]{geometry}

\usepackage{amsmath}
\usepackage{graphicx}
\usepackage[colorlinks=true, allcolors=blue]{hyperref}

\title{OptControl.jl: An Interpreter for Optimal Control Problem}
\author{Jingyi Yang$^1$, Yuebao Yang$^{2}$, Mingtao Li\footnote{Corresponding author}\\Xi'an Jiaotong University}

\begin{document}
\maketitle

\begin{abstract}
OptControl.jl\footnote{https://ai4energy.github.io/OptControl.jl/dev/}(OptControl) implements that modeling optimal control problems with symbolic algebra system based on Julia language, and generates the corresponding numerical optimization codes to solve them with packages from Julia. OptControl does not define a data type, but generates a solution script by handling Julia strings and runs the script automatically. It also provides an interface to save script files. Meanwhile, OptControl supports component-based modeling, which makes it easy to build the optimal control problem of complex systems. All of OptControl's dependency packages come from ecosystem of Julia.
\end{abstract}

\section{Motivation}

The essence of optimal control problem is an optimization problem. More accurately, it is a functional extreme value problem. From the perspective of implementation, the analytical solution of optimal control can be obtained only in certain cases, such as linear systems. However, practical problems usually are nonlinear systems or some much complex systems, and the analytical solutions are generally difficult to obtain. Therefore, the numerical solution that can calculate the result is a sharp weapon. Although the numerical solution will have some deviations in the results, the biased results can also give people some inspiration to understand the real world.

Generally, the numerical solution of optimal control problem can be transformed into numerical optimization problem. Numerical optimization problems can be solved by using JuMP.jl(JuMP). JuMP is a modeling language for mathematical optimization embedded in Julia\cite{DunningHuchetteLubin2017}.

Using JuMP to solve an optimal control problem can be divided into five steps

\begin{itemize}
\item Step 1: Abstract problems from real world,
\item Step 2: Obtain state space model,
\item Step 3: Build numerical optimization model,
\item Step 4: Generate jump model,
\item Step 5: Solve the model.
\end{itemize}

\begin{figure}
\centering
\includegraphics[width=1.0\textwidth]{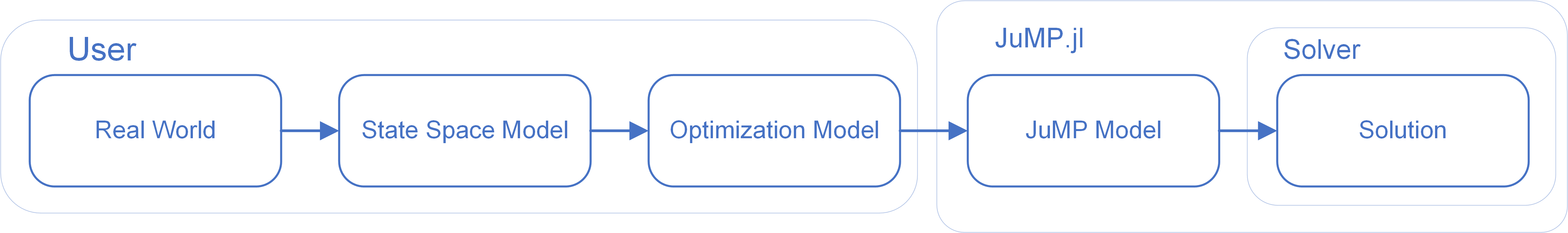}
\caption{\label{fig:f1}Steps to solve an optimal control problem.}
\end{figure}

In fact, JuMP finish Step 4, and the job of Step 5 is finished by a specific solver and JuMP provides an interface to Step 5. Therefore, it can also be considered that JuMP covers Step 4 and Step 5. Then the first three steps are left to the user, as shown in the Figure \ref{fig:f1}.

When dealing with the problem, the first step involves an abstract from real world to mathematical expression, which can only be completed by the advanced human brain. Can steps 2 and 3 be automatically finished? This is exactly what OptControl wants to accomplish. OptControl focuses on automated implementation.

OptControl focuses on automatic implementation, that is, how to automatically build the optimal problem, how to automatically build the jump optimization model and how to call solver to solve automatically. It is not involved that how to solve an optimal problem, how to build a symbolic system and so on .
OptControl aims to integrate resources (various software packages in Julia ecosystem) and completing the above five steps automatically.

So OptControl can be regarded as an interpreter to implement the transformation from state space model to optimal control problem.It has three features:

\begin{enumerate}
\item Accept state space model in form of symbolic system built from Symbolics.jl or ModelingToolkit.jl,
\item Automatically generate the script solving optimization problem expressed in JuMP model and run it automatically,
\item Provide an interface to save script files for users to modify freely
\end{enumerate}

\section{Framework of OptControl.jl}

OptControl's ability is gradually improved.

\subsection{Build optimization model}

\begin{figure}
\centering
\includegraphics[width=1.0\textwidth]{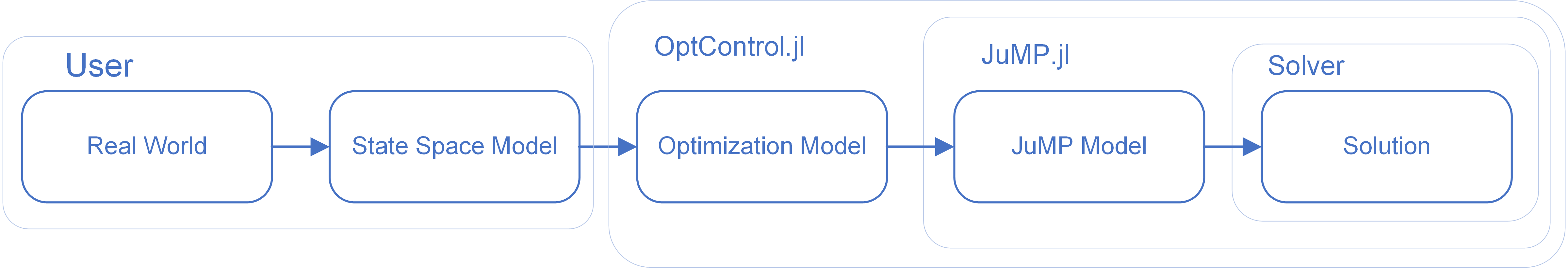}
\caption{\label{fig:f2}Framework of OptControl.jl for solving Step 3 to Step 5.}
\end{figure}

The first thing to complete is the solution from Step 3 to Step 5, as shown in Figure \ref{fig:f2}. The function of this solution is:

\begin{itemize}
\item \textit{generateJuMPcodes}——dealing with linear systems,
\item \textit{generateNLJuMPcodes}——dealing with nonlinear systems.
\end{itemize}

\textit{generateJuMPcodes} and \textit{generateNLJuMPcodes} accept the equation of state expressed in symbolic form with Symbolics. Symbolics.jl(Symbolics) is a symbolic algebra system with high performance and can be extended in user language\cite{gowda2021high}.

The symbolic state space model can be transformed into Julia function which is used to obtain the discrete form of the model. And next taking the state in the discrete model as the variable of jump system. Finally, the JuMP model can be built and problem can be solved by solver.

\subsection{Obtain state space model from ODESystem}

\begin{figure}
\centering
\includegraphics[width=1.0\textwidth]{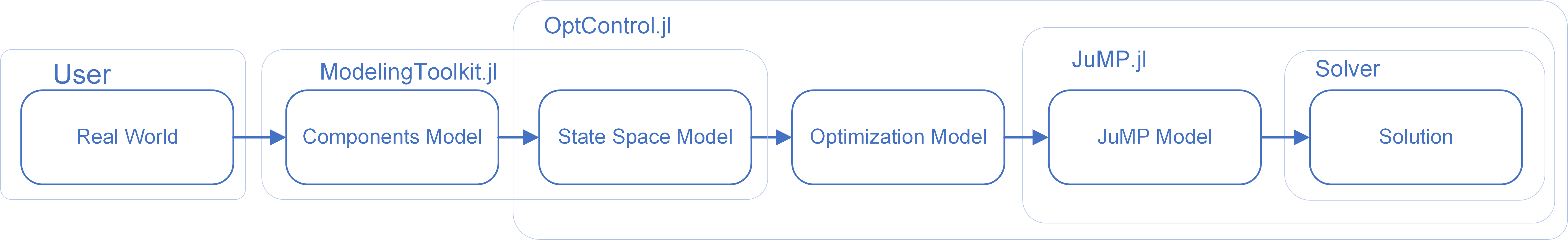}
\caption{\label{fig:f3}Framework of OptControl.jl for solving Step 2 to Step 5.}
\end{figure}

Further, we want to automatically implement Steps 2 to 5 as shown in Figure \ref{fig:f3}. The solutions from step 2 to step 5 need to use the acausal component-based model in ModelingToolkit.jl(ModelingToolkit)\cite{ma2021modelingtoolkit}. ODESystem is a Julia data type in ModelingToolkit for modeling ordinary differential equation system(ODESystem).

The function of this solution is:

\begin{itemize}
\item \textit{generateMTKcodes}—— accept the ODESystem
\end{itemize}

The differential equation governing the system in odesystem is actually the state space form of the optimal control problem. The difference between them is that some of variables in the odesystem system are state-variables in the optimal control problem, while others are control-variables. In other words, the state space equation in the optimal control is essentially the same as the differential equation in ODESystem, but the difference is that some variables have special meanings in optimal control problem.

\textit{generateMTKcodes} uses \textit{generate\_function} in ModelingToolkit to generate Julia function. Then the next step is the same as that in 2.1. Therefore, the solution ideas in 2.1 and 2.2 are the same, because they all generate Julia functions. The generated functions in 2.1 are from Symbolics while the functions in 2.2 are from ODESystem in ModelingToolkit, which is the main difference.

\section{Math work in OptControl}

\subsection{Meaning of Simulation and Control}

The state space equations in the optimal control problem are:

\[\boldsymbol{\dot{x}}(t)=\boldsymbol{A}\boldsymbol{x}(t)+\boldsymbol{B}\boldsymbol{u}(t)=f[\boldsymbol{x}(t),\boldsymbol{u}(t),t]\tag{1}\]

Where $\boldsymbol{x(t)}$ is the state-variable vector of the system, $\boldsymbol{u(t)}$ is the control-variable vector of the system. They are all functions of the independent variable $t$, which means that they change with time. In fact, in the control problem, the coefficient matrix $\boldsymbol{A},\boldsymbol{B}$ can also change with time, expressed as $\boldsymbol{A(t)},\boldsymbol{B(t)}$.

Omitting the $t$, the above equation can be simplified as:

\[\dot{\boldsymbol{x}}=\boldsymbol{A}\boldsymbol{x}+\boldsymbol{B}\boldsymbol{u}=f(\boldsymbol{x},\boldsymbol{u},t)\tag{2}\]

If we consider mathematics without considering physical meaning of $\boldsymbol{x},\boldsymbol{u}$, it is an ordinary differential equation problem with respect to time. If $\boldsymbol{u}$ is not artificially determined, it evolves itself in the system, which will be a problem of solving differential equations, as Equation 3 shows.

\[\dot{\boldsymbol{x}}=f(\boldsymbol{x},\boldsymbol{p},t)\tag{3}\]

Where $\boldsymbol{p}$ is the parameter vector of the system.

Solving the differential equation means dynamic simulation of the system in the real world. Therefore, the essence of control problem and dynamic simulation problem is the same. The difference is whether some variables in the problem can be changed by human intervention. It can also be argued that the dynamic simulation problem is that we want to see how the system evolves, and the control problem is that we want the system to evolve according to our expectations. Because we have expectations, we need human intervention in the system, which is the meaning of $\boldsymbol{u}$. $\boldsymbol{u}$ is a mathematical form of intervention for system.
If we set $\boldsymbol{u}$ without changing it, the problem is still a dynamic simulation problem, because the $\boldsymbol{u}$ does not interfere so that human expectations are not reflected in the system.

The same essence is why OptControl can accept the system from ModelingToolkit. ModelingToolkit was originally a tool for constructing dynamic simulation problems and ODESystem describes the simulation model of a dynamic system, which does not have an interface that can manually intervene in the system - the control variable $\boldsymbol{u}$. We can build an ODESystem and observe how the system changes, without being able to control its evolution direction from beginning to end (in fact, occasional intervention can be achieved through the Callback function, but it is far from the "control"). But we can turn a simulation problem into a control problem with a few changes. Just add meanings to some variables in ODESystem, and the job is done, as the Equation 4 shows.

\[\dot{\boldsymbol{x}}=f(\boldsymbol{x},\boldsymbol{p},t)\Rightarrow \dot{\boldsymbol{x}}=f(\boldsymbol{x},\boldsymbol{u},t)\tag{4}\]

\subsection{Mathematical Form of Optimization}

In the previous section we explored the meaning of $\boldsymbol{u}$. Then the remaining question is how to describe the optimization with math\cite{jorge2006numerical}. The whole optimization problem can be divided into two parts, the optimal in the control process and the optimal in the control final state.Equation 5 represents an objective for the optimal\cite{macki2012introduction}.

\[\min ~~~\Phi(\boldsymbol{x}(t_f),t_f)+\int_{t_0}^{t_f} L[\boldsymbol{x}(t),\boldsymbol{u}(t),t]dt  \tag{5}\]

Where $\Phi(\boldsymbol{x}(t_f),t_f)$ expresses an expectation of the final state. $\int_{t_0}^{t_f} L[\boldsymbol{x}(t),\boldsymbol{u}(t),t]dt$ expresses an expectation in the control process. The combination of Equation 2 and Equation 5 becomes the general form of the governing equation for the optimal control problem.

\[\begin{matrix}
\min~~~~\Phi(\boldsymbol{x}(t_f),t_f)+\int_{t_0}^{t_f} L[\boldsymbol{x}(t),\boldsymbol{u}(t),t]dt\\s.t. \hspace{1.0cm} \dot{\boldsymbol{x}} =
f[\boldsymbol{x}(t),\boldsymbol{u}(t),t] 
\end{matrix} \tag{6}\]

\subsection{Numerical optimization model}

Equation 6 is in continuous form and needs to be discretized if numerical optimization methods are used. The discrete method uses Euler's method\cite{burden2015numerical}, then there are:

\[\begin{matrix}
\min~~~~\Phi(\boldsymbol{x}(t_f),t_f)+\sum_{i=1}^{n} L(\boldsymbol{x}_i,\boldsymbol{u}_i,t_i) \\s.t. \hspace{0.4cm} \boldsymbol{x}_{i+1} =\boldsymbol{x}_{i}+f(\boldsymbol{x}_i,\boldsymbol{u}_i,t_i)*dt
\end{matrix} \tag{7}\]

If the backward Euler's rule is used:

\[\begin{matrix}
\min~~~~\Phi(\boldsymbol{x}(t_f),t_f)+\sum_{i=1}^{n} L(\boldsymbol{x}_i,\boldsymbol{u}_i,t_i) \\s.t. \hspace{0.4cm} \boldsymbol{x}_{i+1} =\boldsymbol{x}_{i}+f(\boldsymbol{x}_{i+1},\boldsymbol{u}_{i+1},t_{i+1})*dt
\end{matrix} \tag{8}\]

In addition, there are many discrete methods, such as trapezoidal method, Adams method\cite{diethelm2004detailed} and so on. Once Equation 6 or Equation 7 is obtained, the next step is to use JuMP to build the corresponding JuMP model, which is the form the solver can solve. These steps are done automatically by OptControl.

\section{Case Study}

\subsection{Case1: Optimal Control of Linear System}

To solve the following linear optimal control problem:

\[\begin{matrix}
 min \int_{0}^{2} u^2dt\\
 s.t. ~~~~~ \dot{\boldsymbol{x}} =\begin{bmatrix}0&1 \\ 0&0\end{bmatrix}\boldsymbol{x}+ \begin{bmatrix}0 \\ 1 \end{bmatrix}u \\
\boldsymbol{x}(0) = \begin{bmatrix} 1 \\ 1 \end{bmatrix}, \boldsymbol{x}(2)=\begin{bmatrix} 0 \\ 0 \end{bmatrix}
\end{matrix}\]

The steps to solve this problem using OptControl are:

\begin{enumerate}
\item Use ModelingToolkit or Symbolics to model the system,
\item Set states such as initial state and final state,
\item Call \textit{generateJuMPcodes} to solve.
\end{enumerate}

\begin{code}
using OptControl, Statistics, ModelingToolkit
@variables t u x[1:2]
f = [0 1; 0 0] * x + [0, 1] * u
L = 0.5 * u^2
t0 = [1.0, 1.0]
tf = [0.0, 0.0]
tspan = (0.0, 2.0)
N = 100
sol = generateJuMPcodes(L, f, x, u, tspan, t0, tf; N=N)
\end{code}

In this optimal problem, the analytical solution of $x_1$ is

$$x_1(t) = 0.5*t^3-1.75*t^2+t+1$$

Comparing the analytical solution and the numerical solution, it can be found the mean square error is 2.696E-6.

Even if the result under this error cannot be used, it is of great reference significance and can give people inspiration. Even if this result is not adopted, it can give people inspiration to solve problem.

\subsection{Case2: Optimal Control of Nonlinear system}

To solve the following nonlinear optimal control problem

\[\begin{array}{c}
min \int_{0}^{2} u^2dt \\ s.t. ~~~~~ \dot{\boldsymbol{x}} =\begin{bmatrix}exp&cos \\ sin&1\end{bmatrix}\boldsymbol{x}+ \begin{bmatrix}0 \\ 1 \end{bmatrix}u \\ \boldsymbol{x}(0) = \begin{bmatrix} 1 \\ 1 \end{bmatrix}, \boldsymbol{x}(1)=\begin{bmatrix} 0 \\ 0 \end{bmatrix}
\end{array}\]
Use ModelingToolkit or Symbolics to define symbolic variables, giving initial and final states. Call \textit{generateNLJuMPcodes} to get the result.

\begin{code}
using OptControl, ModelingToolkit, Test
@variables t u x[1:2]
f = [exp(x[1]) + cos(x[2]), sin(x[1]) + x[2]] + [1, 0] * u
L = u^2
t0 = [1.0, 1.0]
tf = [0.0, 0.0]
tspan = (0.0, 2.0)
N = 100
sol = generateNLJuMPcodes(L, f, x, u, tspan, t0, tf; N=N)
\end{code}

\subsection{Case3: Optimal Control of RC Circuit System}

\begin{figure}
\centering
\includegraphics[width=0.6\textwidth]{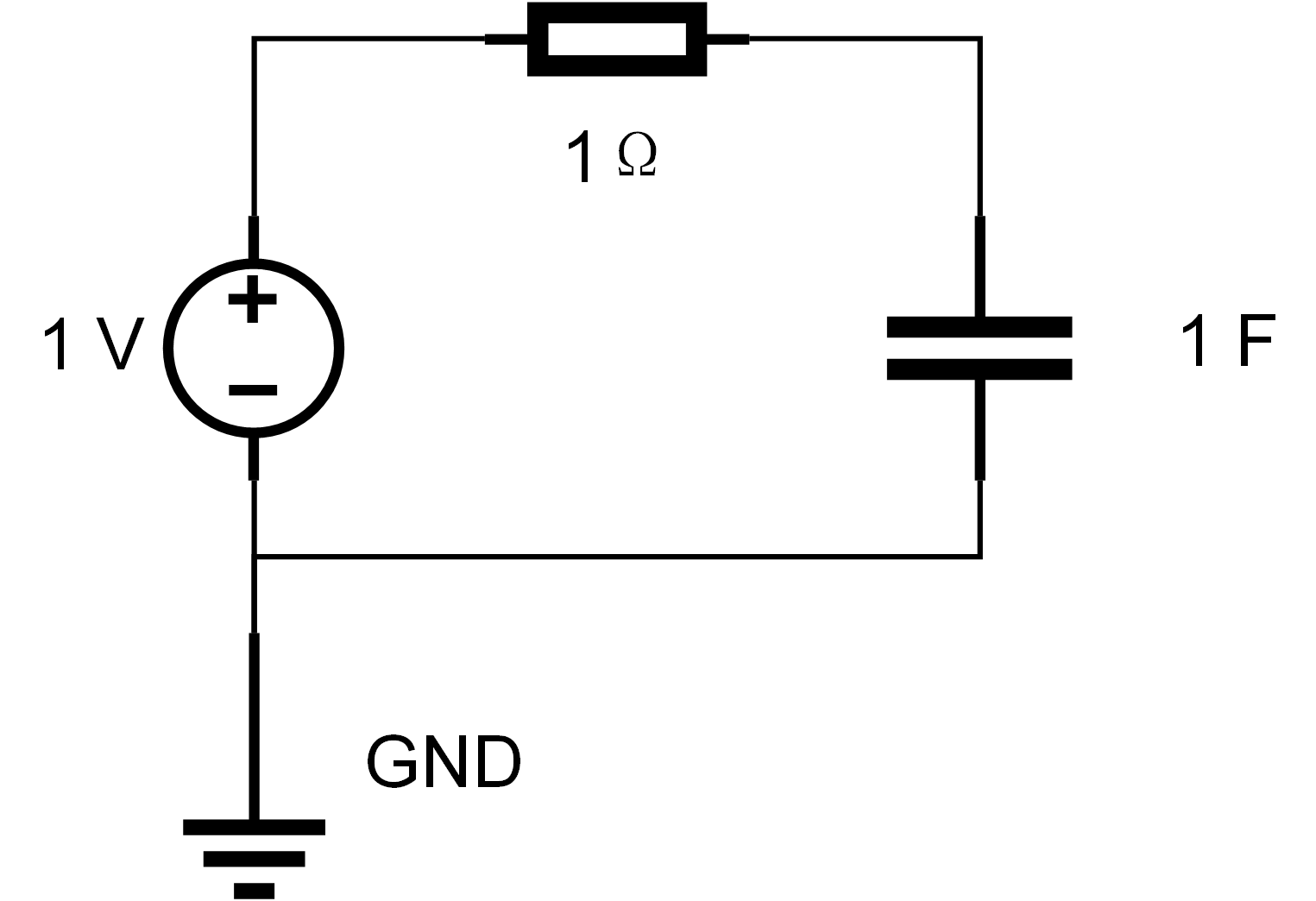}
\caption{\label{fig:f4}RC Circuit System}
\end{figure}

This is a simple circuit system. The supply voltage is 1 volt, the resistance is 1 ohm, and the capacitance is 1 farad, as Figure \ref{fig:f4} shows.

The optimal control problem we constructed is how to change the voltage to make the capacitor voltage change from 1V to 3V within 1s, while satisfying the goal of the voltage as low as possible in the whole process. This has a physical meaning and can illustrate how optimal control can be applied in component-based systems.

The components of the code are taken from the documentation of ModelingToolkit\footnote{https://mtk.sciml.ai/dev/tutorials/spring\_mass/}. If it is a simulation problem, it should be solved by calling DifferentialEquations.jl\cite{rackauckas2017differentialequations} after the ODESystem is simplified. Now it is an optimal control problem, so we set the optimization objective, define the relevant parameters, and call \textit{generateMTKcodes} to solve.

\begin{code}
using OptControl, ModelingToolkit, Test

# Components ......

# Define and Simplify System
R = 1.0
C = 1.0
V = 1.0
@named resistor = Resistor(R=R)
@named capacitor = Capacitor(C=C)
@named source = ConstantVoltage(V=V)
@named ground = Ground()
rc_eqs = [
    connect(source.p, resistor.p)
    connect(resistor.n, capacitor.p)
    connect(capacitor.n, source.n)
    connect(capacitor.n, ground.g)
]
@named _rc_model = ODESystem(rc_eqs, t)
@named rc_model = compose(_rc_model,
    [resistor, capacitor, source, ground])
sys = structural_simplify(rc_model)

# Build Optimal Control Problem and Solve
L = 0.5 * (source.V^2)
t0 = [1.0]
tf = [3.0]
tspan = (0.0, 1.0)
N = 100
sol = OptControl.generateMTKcodes(L, sys, states(sys), [source.V], tspan, t0, tf;N=N)
\end{code}

\section{Conclusion}

OptControl implements the automatic construction from the state space equations to the optimal control problem and the automatic construction from the ordinary differential equation system of the ModelingToolkit to the system optimal control problem. The core of the problem lies in the choice of design control variables $\boldsymbol{u}$. When discretizing, OptControl provides an interface to select the discretization methods. In future work, more discrete methods will be added.

Another important feature of OptControl is that it does not solve the problem directly, but generates a solution script and runs it. This means that OptControl is like a conductor, it decomposes the problem, and then calls the package in the Julia ecosystem to solve the problem. OptControl provides an interface to get the script file, which means that when the script does not meet your needs, you can modify it directly. You can add whatever code you need to the generated script. If you are not familiar with JuMP's modeling language, then you can learn some advanced usage of JuMP through the generated scripts. If you also want to choose some different solvers, modify the script.

In the future, OptControl may provide more interfaces. But it will not change the role as a conductor. OptControl will always focus on automatically generating solutions to optimal control problems, rather than developing a modeling language like ModelingToolkit and JuMP, nor developing solving algorithms like the solvers called by JuMP. The original intention of OptControl is to break through barriers, integrate tools, and solve optimal control problems conveniently and quickly.

\bibliographystyle{unsrt}
\bibliography{sample}

\begin{thebibliography}{1}

\bibitem{DunningHuchetteLubin2017}
Iain Dunning, Joey Huchette, and Miles Lubin.
\newblock Jump: A modeling language for mathematical optimization.
\newblock {\em SIAM Review}, 59(2):295--320, 2017.

\bibitem{gowda2021high}
Shashi Gowda, Yingbo Ma, Alessandro Cheli, Maja Gwozdz, Viral~B Shah, Alan
  Edelman, and Christopher Rackauckas.
\newblock High-performance symbolic-numerics via multiple dispatch.
\newblock {\em arXiv preprint arXiv:2105.03949}, 2021.

\bibitem{ma2021modelingtoolkit}
Yingbo Ma, Shashi Gowda, Ranjan Anantharaman, Chris Laughman, Viral Shah, and
  Chris Rackauckas.
\newblock Modelingtoolkit: A composable graph transformation system for
  equation-based modeling.
\newblock {\em arXiv preprint arXiv:2103.05244}, 2021.

\bibitem{jorge2006numerical}
Nocedal Jorge and J~Wright Stephen.
\newblock Numerical optimization, 2006.

\bibitem{macki2012introduction}
Jack Macki and Aaron Strauss.
\newblock {\em Introduction to optimal control theory}.
\newblock Springer Science \& Business Media, 2012.

\bibitem{burden2015numerical}
Richard~L Burden, J~Douglas Faires, and Annette~M Burden.
\newblock {\em Numerical analysis}.
\newblock Cengage learning, 2015.

\bibitem{diethelm2004detailed}
Kai Diethelm, Neville~J Ford, and Alan~D Freed.
\newblock Detailed error analysis for a fractional adams method.
\newblock {\em Numerical algorithms}, 36(1):31--52, 2004.

\bibitem{rackauckas2017differentialequations}
Christopher Rackauckas and Qing Nie.
\newblock Differentialequations.jl--a performant and feature-rich ecosystem for
  solving differential equations in julia.
\newblock {\em Journal of Open Research Software}, 5(1), 2017.

\end{thebibliography}

\end{document}